\documentclass{article} \usepackage{amsmath,amsthm,amssymb,xspace}
\newcommand{\tiff}{if\textcompwordmark f\xspace}
\newcommand{\abs}[1]{\left|#1\right|}
\newcommand{\gint}[1]{\lfloor#1\rfloor}
\newcommand{\set}[1]{\left\{#1\right\}}
\newcommand{\setst}[2]{\set{\,#1\mid#2\,}}
\newcommand{\func}[3]{#1\colon#2\to#3} \newcommand{\N}{\mathbb{N}}
\newcommand{\Q}{\mathbb{Q}} \newtheorem{prop}{Proposition}
\newtheorem{thm}[prop]{Theorem} \newtheorem{lem}[prop]{Lemma}
\newtheorem{cor}[prop]{Corollary}

\begin{document}
\title{The asymptotic density of dead ends in non-amenable groups}
\author{Andrew~D. Warshall\\Yale University\\Department of
  Mathematics\\P.O. Box 208283\\New Haven, CT
  06520-8283\\\texttt{andrew.warshall@yale.edu}} \date{20 June 2010}
\maketitle

\begin{abstract}
We show that, in non-amenable groups, the density of elements of depth
at least $d$ goes to $0$ exponentially in $d$.
\end{abstract}

Let $G$ be a group, $A$ a finite generating set for $G$. Then the
\emph{depth} (more verbosely the \emph{dead-end depth}) of $g\in G$ is
the distance (in the word metric with respect to $A$) from $g$ to the
complement of the radius-$d_A(1,g)$ closed ball about the
origin. While many common examples of (group, generating set) pairs
admit a uniform bound on the depth of elements, some do not; we say
these groups have \emph{deep pockets}. The standard example of a group
with deep pockets is the lamplighter group with respect to the
standard generating set; see~\cite{CT}.

Although, by definition, groups with deep pockets have elements of
arbitrarily large depth, this leaves open the question of how many
there are. In a remark at p.~91 of \cite{S}, Kyoji Saito asked (in the
context of the study of so-called pre-partition functions) under what
circumstances the asymptotic density of elements of depth $>1$ is
guaranteed to be $0$; he had posed the same question in a personal
communication to the author in 2009.  Inspired by this question (but
as yet unable to answer it), we show that the density of elements of
depth at least $m$ approaches $0$ exponentially in $m$, provided the
group is not amenable.

We begin with a lemma from analysis.

\begin{lem}\label{linml}
Let $\func{f}{(0,1)}{(0,\infty]}$ be completely arbitrary. Then there
  exists $U$ an open interval in $(0,1)$ and $\epsilon>0$ such that,
  for every $x_0\in U$ and $\delta>0$, there is $x\in(0,1)$ with
  $\abs{x-x_0}<\delta$ and $f(x)>\epsilon$.
\end{lem}

\begin{proof}
For every $\epsilon>0$, let $U_\epsilon$ be the set of $x_0$ such
that, for all $\delta>0$, there exists an appropriate $x$. We will
show that some $U_\epsilon$ contains an open interval.

We know that $x_0\in U_{f(x_0)/2}$ for every $x_0$ (just set $x=x_0$),
so $\bigcup_\epsilon U_\epsilon=(0,1)$, since the range of $f$ does
not include $0$.  But, if $\epsilon_1<\epsilon_2$, then
$U_{\epsilon_1}\supseteq U_{\epsilon_2}$, so
\[
\bigcup_\epsilon U_\epsilon=\bigcup_{\epsilon\in\Q}U_\epsilon.
\]
This union is countable, so, since it exhausts $(0,1)$, Baire's
Theorem says that some $U_\epsilon$ is dense on some interval.  But,
for every $\epsilon$, $U_\epsilon$ is closed in $(0,1)$. Thus, if it
is dense on an (open) interval, it includes that interval, so we are
done.
\end{proof}

\begin{cor}\label{linm}
Let $\func{f}{(0,1)}{(0,\infty]}$ be completely arbitrary. Then there
  is $A>1$ such that for every sufficiently large integer $m$ there
  exists $S\subset(0,1)$ such that any two elements of $S$ differ from
  each other, and from $0$ and $1$, by at least $1/m$ and $\prod_{x\in
    S}[1+f(x)]>A^m$.
\end{cor}

\begin{proof}
Choose $U$ and $\epsilon>0$ as in Lemma~\ref{linml}. Let $L$ be the
length of $U$.  If $m\in\N$, then $U$ contains $\gint{Lm}$ disjoint
subintervals of length $1/m$. Let $S$ consist of one point each from
the even-numbered subintervals other than the last (that is from every
other subinterval, omitting the first and last), chosen such that
$f(x)>\epsilon$ for all $x\in S$; this is possible by construction of
$U$.  Then clearly any two elements of $S$ differ from each other, and
from $0$ and $1$, by at least $1/m$. But
\[
\prod_{x\in
  S}[1+f(x)]\ge(1+\epsilon)^{\gint{Lm/2}-1}>(1+\epsilon)^{Lm/2-2}>A^m
\]
for $m$ sufficiently large, so long as $A<(1+\epsilon)^{L/2}$.
\end{proof}

\begin{lem}\label{ddef}
Let $G$ be a group, $A$ a finite generating set for $G$ and $m$,
$n\in\N$.  Let $h\in S_n$.  Then $h$ has depth at least $m$ \tiff $h$
is at distance at least $m$ from $B_n'$ (where by $B_n'$ we mean the
complement of $B_n$ in $G$).
\end{lem}

This is just the definition of depth.

\begin{prop}\label{sphere}
Let $G$ be a non-amenable group and $A$ a finite generating set for
$G$.  Let $B_n$ (respectively $S_n$) be the (closed) ball
(resp.\ sphere) of radius $n$ about the identity in $G$ with respect
to $A$. Let $D_m$ be the set of elements of $G$ of depth at least $m$
with respect to $A$. Then there are $a>0$ and $0<b<1$ such that, for
all $m$, $\sup_{n>m}\abs{S_n\cap D_m}/\abs{B_n}<ab^m$.
\end{prop}

\begin{proof}
It follows from Lemma~\ref{ddef} that, for any $m$,
\[
\sup_{n>m}\frac{\abs{S_n\cap D_m}}{\abs{B_n}}
\]
is the supremum over all $n>m\in\N$ of the fraction of elements of
$B_n$ that are in $S_n$ and at distance at least $m$ from $B_n'$.
Thus there exists a sequence $n_m>m\in\N$ such that the fraction (say
$f_m$) of elements of $B_{n_m}$ that are in $S_{n_m}$ and at distance
at least $m$ from $B_{n_m}'$ is at least $\sup_{n>m}\abs{S_n\cap
  D_m}/\abs{B_n}-\epsilon/2^m$.

I claim that $f_m\le ab^m$ for some $a>0$ and $0<b<1$ independent of
$m$; this will imply the theorem (after increasing $a$ slightly and
possibly $b$).  For each $\delta\in(0,1)$, let $H_{m,\delta}\subset
B_{n_m}$ be the set of points at distance $>1+\delta(m-1)$ from
$B_{n_m}'$.  Every point of $\partial H_{m,\delta}$ is $\in B_{n_m}$
and at distance either $1+\gint{\delta(m-1)}$ or
$2+\gint{\delta(m-1)}$ from $B_{n_m}'$.  For each $\delta$, we know
that, since neither $H_{m,\delta}$ nor any subsequence thereof is F\o
lner,
\[
\liminf_{m\to\infty}\frac{\abs{\partial
    H_{m,\delta}}}{\abs{H_{m,\delta}}}>0,
\]
whence
\[
\inf_m\frac{\abs{\partial
    H_{m,\delta}}}{\abs{H_{m,\delta}}}=l(\delta)>0,
\]
since $n_m>m$ guarantees that $H_{m,\delta}$, hence $\partial
H_{m,\delta}$, is never empty.  Note that $l(\delta)$ is independent
of $m$.

It follows by Corollary~\ref{linm} that there is $A>1$ (also
independent of $m$) such that for all sufficiently large $N\in\N$
there exists $C_N\subset(0,1)$ such that any two elements of $C_N$
differ from each other, and from $0$ and $1$, by at most $1/N$ and
$\prod_{\delta\in C_N}[1+l(\delta)]>A^N$.  We thus have
\begin{itemize}
\item $\prod_{\delta\in C_N}[1+l(\delta)]>A^N$
\item for all $m>2N+1$, $\setst{\partial H_{m,\delta}}{\delta\in C_N}$
  are pairwise disjoint subsets of
  $B_{n_m}-\bigcap_{\delta\in(0,1)}H_{m,\delta}$ and
\item for all $m\in\N$, $\delta\in(0,1)$,
\[
\frac{\abs{\partial H_{m,\delta}}}{\abs{H_{m,\delta}}}\ge l(\delta).
\]
\end{itemize}

It follows that, for all $m>2N+1$, the fraction of elements of
$B_{n_m}$ belonging to
\[
S_{n_m}\cap D_m\subset\bigcap_{\delta\in(0,1)}H_{m,\delta}
\]
is at most
\[
\prod_{\delta\in C_N}\frac{\abs{H_{m,\delta}-\partial
    H_{m,\delta}}}{\abs{H_{m,\delta}\cup\partial
    H_{m,\delta}}}<\prod_{\delta\in
  C_N}\frac{\abs{H_{m,\delta}}}{\abs{H_{m,\delta}}+\abs{\partial
    H_{m,\delta}}}\le\prod_{\delta\in
  C_N}\frac{1}{1+l(\delta)}<A^{-N}.
\]
Since $N$ was arbitrary, we have that, for all $m$, $f_m<A^{(1-m)/2}$,
so the claim is proven and we are done.
\end{proof}

\begin{thm}
Let $G$ be a non-amenable group and $A$ a finite generating set for
$G$.  Let $B_n$ be the (closed) ball of radius $n$ about the identity
in $G$ with respect to $A$. Let $D_m$ be the set of points of depth at
least $m$ with respect to $A$. Then there are $c>0$ and $0<b<1$ such
that, for all $m$,
\[
\lim_{n\to\infty}\frac{\abs{B_n\cap D_m}}{\abs{B_n}}<cb^m.
\]
\end{thm}

\begin{proof}
Since neither the $B_i$ nor any subsequence thereof are F\o lner, we
have
\[
\frac{\abs{B_{i-2}}}{\abs{B_i}}<L<1
\]
for all $i$, where $L$ depends only on $G$ and $A$. Thus
$\abs{B_{i-j}}/\abs{B_i}<L^{\gint{j/2}}<L^{j/2-1}$, since certainly
$\abs{B_i}\le\abs{B_{i+1}}$ for all $i$. Thus, for any $m<n$,
\begin{multline*}
\frac{\abs{B_n\cap D_m}}{\abs{B_n}}=\sum_{i=0}^n\frac{\abs{S_i\cap
    D_m}}{\abs{B_n}}=\sum_{i=0}^n\frac{\abs{S_i\cap
    D_m}}{\abs{B_i}}\frac{\abs{B_i}}{\abs{B_n}}\\=\frac{\abs{B_m}}{\abs{B_n}}+\sum_{i=m+1}^n\frac{\abs{S_i\cap
    D_m}}{\abs{B_i}}\frac{\abs{B_i}}{\abs{B_n}}<\frac{\abs{B_m}}{\abs{B_n}}+\sum_{i=m+1}^nab^m\frac{\abs{B_i}}{\abs{B_n}}\\<\frac{\abs{B_m}}{\abs{B_n}}+\sum_{i=m+1}^nab^mL^{(n-i)/2-1}<\frac{\abs{B_m}}{\abs{B_n}}+ab^m\sum_{i=0}^\infty
L^{i/2-1}=\frac{\abs{B_m}}{\abs{B_n}}+cb^m/2,
\end{multline*}
where the first inequality is by Proposition~\ref{sphere} and the
second by the preceding sentence. (Here $c=2a\sum_{i=0}^\infty
L^{i/2-1}$ depends only on $G$ and $A$.) Thus
\[
\limsup_{n\to\infty}\frac{\abs{B_n\cap D_m}}{\abs{B_n}}\le cb^m/2<cb^m
\]
for all $m$, as claimed.
\end{proof}

\end{document}